\documentclass[12pt]{amsart}
\usepackage{amssymb,slashed, color,array}
\usepackage{amsmath,amsfonts,amsthm,euscript,mathrsfs,multirow}
\usepackage[all]{xy}
\usepackage[english]{babel}
\usepackage{epsf}
\usepackage{graphicx}
\csname @addtoreset\endcsname{equation}{section}

\newtheorem{theorem}{Theorem}[section]
\newtheorem{lemma}[theorem]{Lemma}

\newtheorem{prop}[theorem]{Proposition}

\theoremstyle{definition}

\theoremstyle{remark}

\newenvironment{notation and conventions}{\textbf{Notation and conventions.}}{ }
\DeclareFontFamily{U}{rsf}{} \DeclareFontShape{U}{rsf}{m}{n}{ <5> <6> rsfs5 <7> <8> <9> rsfs7 <10-> rsfs10}{}
\DeclareMathAlphabet\Scr{U}{rsf}{m}{n}
\definecolor{pink}{rgb}{1,0,1}

  \usepackage[pdftex]{hyperref}

\begin{document}

\begin{titlepage}
\begin{center}
\baselineskip=14pt{\LARGE
 On stringy invariants of GUT vacua\\
}
\vspace{1.5cm}
{\large  James Fullwood$^{\spadesuit,}$ , Mark van Hoeij$^\clubsuit$\footnote{Supported by NSF Grant number 1017880.}
  } \\
\vspace{.6 cm}

${}^\spadesuit$Department of Mathematics, The University of Hong Kong, Pok Fu Lam Road, Hong Kong, China\\
${}^\clubsuit$Department of Mathematics, Florida State University, Tallahassee, FL 32306, U.S.A.\\

\end{center}

\vspace{1cm}
\begin{center}

{\bf Abstract}
\vspace{.3 cm}
\end{center}

{\small
We investigate aspects of certain stringy invariants of singular elliptic fibrations which arise in engineering Grand Unified Theories in F-theory. In particular, we exploit the small resolutions of the total space of these fibrations provided recently in the physics literature to compute `stringy characteristic classes', and find that numerical invariants obtained by integrating such characteristic classes are predetermined by the topology of the base of the elliptic fibration. Moreover, we derive a simple (dimension independent) formula for pushing forward powers of the exceptional divisor of a blowup, which one may use to reduce any integral (in the sense of Chow cohomology) on a small resolution of a singular elliptic fibration to an integral on the base. We conclude with a speculatory note on the cohomology of small resolutions of GUT vacua, where we conjecture that certain simple formulas for their Hodge numbers may be given solely in terms of the first Chern class and Hodge numbers of the base.}

\addtocounter{page}{1}

\newpage

 \tableofcontents{}

\vfill

\noindent $^\spadesuit$Email: {fullwood at maths.hku.hk}\\
$^\clubsuit$Email: {hoeij at math.fsu.edu}\\

\end{titlepage}

\section{Introduction}\label{intro}
F-theory provides a geometric platform for engineering Grand Unified Theories within the framework of string theory \cite{Vafa}\cite{BIKMSV}\cite{KatzVafa}. The geometric apparatus of F-theory is an elliptic fibration $\varphi:Y\to B$, whose total space $Y$ is a Calabi-Yau fourfold and whose base $B$ is a compact smooth algebraic variety over $\mathbb{C}$. Crucial to the geometry and the associated physics of the fibration are its singular fibers, which lie over a hypersurface in the base referred to as the \emph{discriminant locus} of the fibration. Physicallly significant aspects of the singular fibers include the fact that they encode the structure of gauge theories associated with D-branes wrapping components of the discriminant locus over which they appear. The standard procedure for the realization of a desired gauge group/singular fiber is to introduce a Weierstrass elliptic fibration in Tate form:

\[
Y:(y^2z+a_{1}xyz+a_3yz^2=x^3+a_2x^2z+a_4xz^2+a_6z^3)\subset \mathbb{P}(\mathscr{E}),
\]
\\
where $\pi:\mathbb{P}(\mathscr{E})\to B $ is the projective bundle of \emph{lines} in $\mathscr{E}=\mathscr{O}_B\oplus \mathscr{L}^2 \oplus \mathscr{L}^3$ for some suitably ample line bundle $\mathscr{L}\to B$ \footnote{For $Y$ to be Calabi-Yau one must take $\mathscr{L}$ to be the anti-canonical bundle $\mathscr{O}(-K)\to B$.}. Choosing $x$, $y$, and $z$ to be respective sections of $\mathscr{O}_{\mathbb{P}(\mathscr{E})}(1)\otimes \pi^*\mathscr{L}^2$, $\mathscr{O}_{\mathbb{P}(\mathscr{E})}(1)\otimes \pi^*\mathscr{L}^3$, $\mathscr{O}_{\mathbb{P}(\mathscr{E})}(1)$, and each $a_i$ to be a section of $\pi^*\mathscr{L}^i$ then realizes $Y$ as the zero-scheme associated with the vanishing of a section of $\mathscr{O}_{\mathbb{P}(\mathscr{E})}(3)\otimes \pi^*\mathscr{L}^6$, and naturally determines a proper surjective morphism $\varphi:Y\to B$ such that the generic fiber is an elliptic curve. When the $a_i$s are suitably generic $Y$ is a smooth hypersurface in $\mathbb{P}(\mathscr{E})$ with nodal and cuspidal cubics as its singular fibers, which correspond to abelian gauge groups in the physical theory. Such an elliptic fibration we will refer to as a \emph{smooth Weierstrass fibration}. To realize non-abelian gauge groups, one resorts to Tate's algorithm \cite{Tate}, which renders conditions on the coefficient sections $a_i$ which inflict singularities on the total space of the fibration $Y$ in such a way that a resolution of the codimension one singularities of $Y$ produce the desired singular fiber over a certain component of the discriminant locus $\Delta$, which is given by

\[
\Delta:(4F^3+27G^2=0)\subset B,
\]
\\
where

\[
\begin{cases} F=-\frac{1}{48}(b_2^2-24b_4)\\ G=-\frac{1}{864}(36b_2b_4-b_2^3-216b_6) \\
b_2=a_1^2+4a_2 \\ b_4=a_1a_3+2a_4 \\ b_6=a_3^2+4a_6. \end{cases}
\]
\\
For example, to engineer an $\text{SU}(n)$ gauge group one imposes that the $a_i$s vanish to certain orders along a divisor $D_{\text{GUT}}:(w=0)\subset B$ in such a way that the total space $Y$ is inflicted with a surface of $A_{n-1}$ singularities over $D_{\text{GUT}}$. The discriminant locus then factors as

\[
\Delta: (w^n\cdot g=0)\subset B,
\]
\\
where $\Delta':(g=0)$ is generically irreducible. Such an elliptic fibration we will refer to as an \emph{$\text{SU}(n)$ elliptic fibration}. Tate's algorithm then ensures that a split $I_n$ fiber (a chain of $n$ rational curves) appears generically over the divisor $D_{\text{SU}(n)}$ upon a resolution of the codimension one singularities of $Y$. In this note we consider the cases of $\text{SU}(5)$, $\text{SO}(10)$ and $\text{E}_6$ vacua. For the $\text{SO}(10)$ and $\text{E}_6$ cases, the coefficients of the Tate form are tweaked in such a way that a split $I_1^*$ fiber (in the SO(10) case) and split $IV^*$ fiber (in the $\text{E}_6$ case)  appears over $D_{\text{GUT}}$ after a resolution of the codimension one singularities of $Y$, as is necessary for the associated physical theory. But due to the exotic nature of singularities, a crepant resolution of singularities in \emph{all} codimensions must transpire for a well-defined physical theory to be associated with the fibration. Thus enters the theory of singular varieties, their resolutions, and associated invariants. When the base of the fibration is a toric variety, so too is the ambient projective bundle $\mathbb{P}(\mathscr{O}_B\oplus \mathscr{L}^2 \oplus \mathscr{L}^3)$ in which it resides, thus resolution procedures via toric methods are readily available \cite{CandelasFont}\cite{BGJW}. Over non-toric bases the algorithmic methods of toric geometry no longer apply, and one must perform an honest resolution by hand. As it is not clear at present which bases give rise to phenomenologically realistic vacua, we prefer to work in the general setting where the only assumption we make on the base is that such an elliptic fibration exists. In the case of $\text{SU}(5)$ models Esole and Yau performed an explicit resolution procedure yielding six distinct small (and thus crepant) resolutions of the total space of the fibration which differ by flop transitions \cite{EY}. Most notably, they make no assumption on the base and do not assume any Calabi-Yau hypothesis. They also provided a detailed analysis of the singular fiber structure of the resolution varieties, and found that the structure of enhancement was different than what had been previously conjectured in the physics literature (and found fibers not on the list of Kodaira as well). For SO(10) and $\text{E}_6$ fibrations, small resolution procedures over arbitrary bases were recently presented in the physics literature \cite{SO(10)}\cite{E6}. 

Invariants associated with a small resolution of a `GUT' singular elliptic fibration may be viewed as invariants associated with the original \emph{singular} fourfold $Y$, following the convolution of mirror symmetry and algebraic geometry which has congealed into a theory of `stringy invariants' of singular varieties. From this perspective a resolution variety is an auxiliary object which allows us to glean information from the original, singular variety. In particular, if $f:\widetilde{X}\to X$ is a crepant resolution\footnote{That is, $f^*K_X=K_{\widetilde{X}}$} of a (not too) singular variety $X$ then invariants of $\widetilde{X}$ are often deemed `stringy invariants' of $X$ provided the invariant is independent of the chosen crepant resolution. For example, Kontsevich showed using motivic integration that the Hodge numbers $h^{p,q}$ of a crepant resolution are independent of the chosen crepant resolution \cite{Kontsevich}, leading Batyrev to define a notion of \emph{stringy Hodge numbers} \cite{Batyrev}. But since crepant resolutions don't always exist, intrinsic definitions are sought of stringy invariants which agree with the associated invariant of a crepant resolution when one exists. Somewhat recently the notion of \emph{stringy Chern classes} were independently defined by Aluffi and de Fernex \emph{et al} \cite{Aluffi}\cite{deFernex}. Aluffi's approach uses his theory of `celestial integration' to define stringy Chern classes while de Fernex \emph{et al} use Kontsevich's theory of motivic integration, yet the seemingly two different approaches reproduce the same class\footnote{Both approaches require that we work over a field of characteristic zero.}. An indication that the moniker `stringy Chern class' is justified lies in the fact that integrating $c_{str}(X)$ yields the \emph{stringy Euler characteristic} $\chi_{str}(X)$ as defined by Batyrev \cite{Batyrev2}, i.e.,

\[
\int_{X}c_{str}(X)=\chi_{str}(X),
\]
\\
which we may view as a `stringy' version of the Poincar\'{e}-Hopf theorem. Moreover, if $f:\widetilde{X}\to X$ is a crepant resolution then $c_{str}(X)=f_*c(\widetilde{X})$ ($f_*$ is the proper push forward map associated with $f$) and $c_{str}(X)=c(X)$ for smooth $X$, as should any appropriate generalization of Chern class. If $X$ is a closed subvariety of some smooth variety $M$ and $N_{X}M$ denotes the normal cone of $X$ in $M$, then the class

\[
\frac{c(TM)}{c(N_{X}M)}\cap [X]
\]
\\
is denoted by $c_{\text{F}}(X)$ and is referred to as the \emph{Fulton class} of $X$, and this class also coincides with $c(X)$ when $X$ is smooth. The difference $c_{str}(X)-c_{\text{F}}(X)$ we will refer to as the \emph{stirngy Milnor class} of $X$, which we will denote by $\mathcal{M}_{str}(X)$ (compare with the definition of Milnor class \cite{Aluffi2}). The stringy Milnor class is then a `stringy' invariant supported on the singular locus of $X$ which measures the deviation of $c_{str}(X)$ from the total homology Chern class of a smooth variety in the same rational equivalence class as $X$. For one interested in numerical invariants then $\int_{X} \mathcal{M}_{str}(X)$ measures the deviation of $\chi_{str}(X)$ from the topological Euler characteristic of a smooth variety in the same rational equivalence class as $X$.

For $\varphi:Y\to B$ a GUT elliptic fibration the integer $\int_{Y} \mathcal{M}_{str}(Y)$ is physically relevant for the computation of 3-brane tadpoles, as it modifies the Euler characteristic of the total space of a smooth Weierstrass fibration $\psi:Z\to B$ to obtain the Euler characteristic of a crepant resolution of $Y$. In the case of SU(5), SO(10) and $\text{E}_6$ fibrations, we relate this invariant of $Y$ to invariants of the base $B$ via the proper pushforward map $\varphi_*$, as $\int_{Y} \mathcal{M}_{str}(Y)=\int_{B}\varphi_* \mathcal{M}_{str}(Y)$. We recall that for $\psi:Z\to B$ a smooth Weierstrass fibration it is known that

\[
\psi_*c(Z)=\frac{12L}{1+6L}c(B),
\]
\\
where $L$ is the first Chern class of the line bundle $\mathscr{L}\to B$ used to define a smooth Weierstrass fibration \cite{AE1}. As such, we may arrive at an expression for $\varphi_*\mathcal{M}_{str}(Y)$ by subtracting $\frac{12L}{1+6L}c(B)$ from $\varphi_*c_{str}(Y)$, the outcome of which we record in the following

\begin{prop}\label{SVW} Let $\varphi:Y\to B$ be an \emph{SU(5)}, \emph{SO(10)} or $\text{\emph{E}}_6$ elliptic fibration with $\text{\emph{dim}}(B)\leq3$, $D=D_{\emph{GUT}}\in A^*B$ and let $L=c_1(\mathscr{L})$. Then $\varphi_*\mathcal{M}_{str}(Y)$ is a multiple (in $A^*B$) of $c(D)$, in particular,

\[
\tag{*}
\varphi_*\mathcal{M}_{str}(Y)=W\cdot c(D), 
\]

\vspace{0.15in}

\noindent where 
\[
W=\begin{cases} \frac{5(36L^3+42L^2+16L-31LD-30L^2D-6D)}{(1+6L)(1+6L-5D)(1+L)} \hspace{1.33in}  \emph{(SU(5) case)}  \\  \frac{4(108L^3+84L^2+21L+10D^2+45D^2L-77DL-144DL^2-8D)}{(1+6L)(1+6L-5D)(1+2L-D)} \hspace{0.37in} \emph{(SO(10) case)} \\ \frac{3(252L^3+162L^2-378L^2D-165LD+30L-12D+140LD^2+30D^2)}{(1+6L)(1+6L-5D)(1+3L-2D)} \hspace{0.32in} \emph{(E}_6 \emph{ case).} \end{cases}
\]
\\
\end{prop}

In the SU(5) and $\text{E}_6$ cases, $\int_{B} \varphi_*\mathcal{M}_{str}(Y)$ over a base of dimension three coincides with $\chi_{str}(Y)-\chi(Z)$ for $Z$ the total space of a smooth Weierstrass fibration as computed in the physics literature \cite{MS}\cite{E6}. We also note that the assumption on the dimension of the base is required as the small resolutions we exploit to compute stringy Chern classes resolve the singularities of $Y$ only up to codimension three in the base. 

From the F-theory perspective, many physically significant quantities may be expressed in terms of integrals (or rather, intersection numbers) on a crepant resolution $\widetilde{Y}$ of the singular fourfold $Y$. In each of the small resolutions we consider, a resolution variety $\widetilde{Y}$ results from the taking the proper transform of the original singular elliptic fibration $Y$ under a sequence of at least four blowups, thus keeping track of intersection data on $\widetilde{Y}$ can be rather unpleasant \cite{MS}\cite{E6}. As such, the methods used in deriving Proposition~\ref{SVW} (namely Lemmas \ref{BC}, \ref{BC2} and the push forward formula of \cite{James}) may be used to bypass intersection theory on $\widetilde{Y}$ and reduce \emph{any} integral (in the sense of Chow cohomology) on $\widetilde{Y}$ to an integral on the base $B$. In what follows we do not assume the total space of the fibration is Calabi-Yau, as the Calabi-Yau case is easily recovered by letting $\mathscr{L}=\mathscr{O}(-K_B)$.
\\

\emph{Acknowledgements.} JF would like to thank Mboyo Esole not only for the motivation to initiate this project, but for the many useful discussions and insights shared concerning the geometry of elliptic fibrations. JF would also like to thank Paolo Aluffi for his constant support and influence. MvH  was supported by NSF Grant number 1017880 during the course of this project.

\section{A little blowup calculus}\label{blowup calculus} In this section we recall a lemma of Aluffi regarding Chern classes of blowups along with a derivation of a formula for pushing forward powers of the exceptional divisor of a blowup under (somewhat) mild assumptions. We combine both results in the following

\begin{lemma}[Aluffi,\hspace{0.07in}$\cdot$\hspace{0.07in}]\label{BC} Let $f:\widetilde{X}\to X$ be the blowup of a smooth variety $X$ along a smooth complete intersection $V:(F_1=F_2=\cdots =F_k=0)\subset X$, let $E\in A^*\widetilde{X}$ be the class of the exceptional divisor and let $U_i\in A^{*}X$ be the class of $F_i=0$. Then

\[
\tag{\dag}
c(T_{\widetilde{X}})=\frac{(1+E)(1+f^*U_1-E)\cdots (1+f^*U_k-E)}{(1+f^*U_1)\cdots (1+f^*U_k)}\cdot f^*c(T_X),
\]
\\
and

\[
\tag{\dag\dag}
f_*(E^n)=\sum_{i=1}^k\left(\prod_{j\neq i}\frac{U_j}{U_j-U_i}\right)U_i^n
\]
\\
for all $n\geq0$.
\end{lemma}

As $U_j-U_i$ is not necessarily invertible in $A^*X$, the quantities $\frac{1}{U_j-U_i}$ appearing in formula ($\dag\dag$) are formal objects which end up (formally) canceling to give a well-defined class in $A^*X$. We prove only formula ($\dag\dag$) as a proof of ($\dag$) can be found in \cite{Aluffi3}. The proof is merely a simple observation about the \emph{Segre class} of $V$ in $X$, which is denoted $s(V,X)$ (for more on Segre classes, see \cite{Fulton}).

\begin{proof} By the birational invariance of Segre classes $f_*s(E,\widetilde{X})=s(V,X)$, where $f_*$ is the proper pushforward associated with the map $f:\widetilde{X}\to X$. Then from the fact that $E$ and $V$ are both regularly imbedded in $\widetilde{X}$ and $X$ respectively, we have

\[
s(E,\widetilde{X})=\frac{E}{1+E}=E-E^2+E^3-\cdots
\]
and
\[
s(V,X)=\prod_{i=1}^k\frac{U_i}{1+U_i},
\]
thus
\[
f_*(E-E^2+E^3-\cdots)=\prod_{i=1}^k\frac{U_i}{1+U_i}.
\]
\\
Matching terms of like dimension we see that $f_*(E^n)$ equals the coefficient of $t^n$ in the series $(-1)^{k+1}\prod_{i=1}^k\frac{tU_i}{1+tU_i}$, which (one can prove by induction) is precisely $\sum_{i=1}^k\left(\prod_{j\neq i}\frac{U_j}{U_j-U_i}\right)U_i^n$.
\end{proof}

A nice feature of formula ($\dag\dag$) is that it gives a dimension independent way of pushing forward classes in the Chow ring of the blowup $A^*\widetilde{X}$, i.e., a class that's given in terms of a rational expression such as formula $(\dag)$ may be pushed forward in terms of another rational expression. For example, assume we are under the hypotheses of Lemma~\ref{BC}, let $k=2$, $U_1=U$, $U_2=V$ and let $T\in A^*(X)$ \footnote{Here and throughout, we often fail to distinguish between classes and their pullbacks.}. Then for $n\geq0$ we have

\begin{eqnarray*}
f_*\left(\frac{E^n}{1+T-E}\right)&=&f_*\left(\frac{E^n}{1+T}\left(\frac{1}{1-\frac{E}{1+T}}\right)\right) \\
                    &=&f_*\left((1+T)^{n-1}\displaystyle\sum_{m=n}^{\infty}\left(\frac{E}{1+T}\right)^m\right) \\
                    &\overset{\text{by}\hspace{0.05in} (\dag\dag)}=&\frac{(1+T)^{n-1}}{V-U}\left(V\displaystyle\sum_{m=n}^{\infty}\left(\frac{U}{1+T}\right)^m-
                    U\displaystyle\sum_{m=n}^{\infty}\left(\frac{V}{1+T}\right)^m\right) \\
                    &=&\frac{(1+T)^{n-1}}{V-U}\left(\left(\frac{U}{1+T}\right)^n\frac{V}{1-\frac{U}{1+T}}-\left(\frac{V}{1+T}\right)^n\frac{U}{1-\frac{V}{1+T}}\right) \\
                    &=&\frac{1}{V-U}\left(\frac{U^nV}{1+T-U}-\frac{V^nU}{1+T-V}\right).
\end{eqnarray*}
\\
A general formula for arbitrary $k$ is derived in a similar fashion, which we record in the following

\begin{lemma}\label{BC2} Under the assumptions of Lemma~\ref{BC}, let $T\in A^*(X)$. Then

\[
f_*\left(\frac{E^n}{1+T\pm E}\right)=\sum_{i=1}^k\left(\prod_{j \neq i} \frac{U_j}{U_j - U_i}\right)\frac{U_i^n}{1+T\pm U_i}
\]
\\
for $n\geq0$.
\end{lemma}

As one can see by glancing at formula ($\dag$), such pushforward formulas are all that is needed to pushforward Chern classes of blowups, whose derivation involves only simple manipulations of rational expressions and geometric series. Moreover, such manipulations are manifestly independent of the dimension of $X$.

The small resolutions of GUT vacua we exploit to compute stringy characterstic classes are all obtained by a sequence of blowups along smooth complete intersections, and then taking the proper transform of the total space of the singular elliptic fibration $Y$ under the blowups. As such, each blowup in the resolution satisfies the hypotheses of Lemma~\ref{BC}. Thus if $\widetilde{Y}$ is a small resolution of $Y$ then pushing forward a class $\gamma\in A^*\widetilde{Y}$ to $A^*B$ amounts to applying of formula ($\dag\dag$) until one arrives at a class in $A^*\mathbb{P}(\mathscr{E})$, and then applying the pushforward formula for classes in a projective bundle which was first derived in \cite{James}.

\section{The fibrations under consideration}\label{fib} Let $B$ be a smooth compact complex algebraic variety of dimension at most three endowed with a line bundle $\mathscr{L}\to B$. As mentioned in \S\ref{intro}, Grand Unified Theories are engineered in F-theory via an elliptic fibration in Tate form:

\[
Y:(y^2z+a_{1}xyz+a_3yz^2=x^3+a_2x^2z+a_4xz^2+a_6z^3)\subset \mathbb{P}(\mathscr{E}),
\]
\\
where $\pi:\mathbb{P}(\mathscr{E})\to B $ is the projective bundle of \emph{lines} in $\mathscr{E}=\mathscr{O}_B\oplus \mathscr{L}^2 \oplus \mathscr{L}^3$, the $a_i$s are global sections of $\mathscr{L}^i$ and $x$, $y$, and $z$ are chosen such that the equation for $Y$ corresponds to the zero-locus of a section of $\mathscr{O}_{\mathbb{P}(\mathscr{E})}(3)\otimes \pi^*\mathscr{L}^6$. The physics of F-theory requires that $Y$ be an anti-canonical divisor of $\mathbb{P}(\mathscr{E})$ (which is achieved by taking $\mathscr{L}=\mathscr{O}(-K_B)$), but is not necessary for our considerations thus we make no assumption on $\mathscr{L}$ other than the fact that suitably generic $a_i$s exist. Such a hypersurface naturally determines a proper surjective morphism $\varphi:Y\to B$ such that the generic fiber is an elliptic curve. The elliptic fiber degenerates to a nodal cubic over a generic point of the the discriminant locus $\Delta\subset B$, and enhances to a cusp over a curve in $\Delta$ \footnote{We recall that the precise definition of $\Delta$ was given in \S\ref{intro}.}. As nodal and cuspidal singular fibers correspond to abelian gauge groups in F-theory, one engineers a particular non-abelian gauge group by introducing singularities into the total space $Y$ of the elliptic fibration in such a way that a resolution of the codimension one singularities of $Y$ will result in the corresponding singular fiber appearing over a divisor $D_{\text{SU}(5)}:(w=0)$ of the base. In particular, Tate's algorithm prescribes the precise orders of vanishing along $D_{\text{GUT}}$ that each coefficient section $a_i$ must satisfy in order to realize a particular gauge group. For example in the SU(5) case, Tate's algorithm renders the following (re)definitions of the $a_i$s (for $i\neq 1$):

\[
a_2=\beta_4w, \quad a_3=\beta_3w^2, \quad a_4=\beta_2w^3, \quad a_6=\beta_0w^5.
\]
\\
Each $\beta_j$ is then necessarily a section of $\mathscr{L}^{6-j}\otimes \mathscr{L}_{\text{GUT}}^{j-5}$ ($\mathscr{L}_{\text{GUT}}$ is the line bundle corresponding to the divisor $D_{\text{GUT}}$) and the new equation for $Y$ then becomes

\[
Y:(y^2z+a_1xyz+\beta_3w^2yz^2=x^3+\beta_4wx^2z+\beta_2w^3xz^2+\beta_0w^5z^3)\subset \mathbb{P}(\mathscr{E}).
\]
\\
Upon such redefinitions of the $a_i$s we will refer to such a fibration as an \emph{$\text{SU}(5)$ elliptic fibration}. 

For the SO(10) case, Tate's algorithm prescribes the same definitions of the $a_i$s except $a_1$, which is now required to vanish to order one along $D_{\text{GUT}}$, i.e., $a_1=\vartheta w$, where $\vartheta$ is a generic section of $\mathscr{L}\otimes \mathscr{L}_{\text{GUT}}^{-1}$ independent from $\beta_4$. The new equation for $Y$ then becomes

\[
Y:(y^2z+\vartheta wxyz+\beta_3w^2yz^2=x^3+\beta_4wx^2z+\beta_2w^3xz^2+\beta_0w^5z^3)\subset \mathbb{P}(\mathscr{E}).
\]
\\
We will refer to such a fibration as an \emph{SO(10) elliptic fibration}.

For the $\text{E}_6$ case, Tate's algorithm prescribes the same definitions of the $a_i$s as in the SO(10) case except $a_2$, which is now required to vanish to order two along $D_{\text{GUT}}$, i.e. $a_2=\eta w^2$, where $\eta$ is a generic section of $\mathscr{L}^2\otimes \mathscr{L}_{\text{GUT}}^{-2}$ independent from $\beta_3$. The new equation for $Y$ then becomes

\[
Y:(y^2z+\vartheta wxyz+\beta_3w^2yz^2=x^3+\eta w^2x^2z+\beta_2w^3xz^2+\beta_0w^5z^3)\subset \mathbb{P}(\mathscr{E}).
\]
\\
We will refer to such a fibration as an \emph{$E_6$ elliptic fibration}\footnote{We apologize for the possible confusion, as these `$\text{E}_6$ elliptic fibrations' are different from the smooth family of elliptic fibrations previously referred to in the physics literature also as `$\text{E}_6$ elliptic fibrations' \cite{KLRY}\cite{AE2}.}.

The singular loci of $\text{SU}(5)$, SO(10) and $\text{E}_6$ elliptic fibrations coincide with the smooth complete intersection\footnote{Though the singular loci of SU(5), SO(10) and $\text{E}_6$ elliptic fibrations coincide, it is the \emph{scheme} structure on the singular locus of each fibration determined by the ideal generated by the partial derivatives of a local defining equation which distinguishes the singularities from one another.}

\[
Y_{sing}:(x=y=w=0)\subset \mathbb{P}(\mathscr{E}).
\]
\\

\section{Stringy Chern classes}
As mentioned in \S\ref{intro}, stringy Chern classes were defined independently by Aluffi and de Fernex \emph{et al} using two different technologies which produce the same class. Aluffi's approach was in terms of his theory of `celestial integration', while de Fernex \emph{et al} make use of Kontsevich's theory of motivic integration\footnote{As Aluffi often does, we put quotes around celestial integration as it is not defined with respect to some sort of measure, and thus is not an honest `integral'. However, it has many points of contact with motivic integration and has integral-like properties such as a change of variable formula with respect to birational maps.}. When a crepant resolution $f:\widetilde{X}\to X$ of a singular variety $X$ exists, we have $c_{str}(X)=f_*c(\widetilde{X})$, and is independent of the resolution. As such, we may now exploit the resolution procedures recently provided in the physics literature (which we outline in the Appendix) along with Lemmas 2.1 and 2.2 to compute stringy Chern classes of SU(5), SO(10) and $\text{E}_6$ elliptic fibrations, then we push forward these classes to the base using the pushforward formula of \cite{James}. More precisely, for $\varphi:Y\to B$ an SU(5), SO(10) or $\text{E}_6$ elliptic fibration and $f:\widetilde{Y}\to Y$ a small resolution of $Y$ we compute $\widetilde{\varphi}_*c(\widetilde{Y})$, where $\widetilde{\varphi}=\varphi \circ f$. We give details of the SU(5) case as the SO(10) and $\text{E}_6$ cases follow \emph{mutatis mutandis}.

Let $\widetilde{Y}$ be one of the small resolutions of an SU(5) elliptic fibration as given in \S A.1, $H=c_1(\mathscr{O}_{\mathbb{P}(\mathscr{E})}(1))$ and let $L=c_1(\mathscr{L})$. By adjuntion and Lemma~\ref{BC} ($\dag$) we have

\begin{eqnarray*}
\iota_*c(\widetilde{Y})&=&\frac{(1+E_4)(1+\mathcal{Y}_2-E_4)(1+V-E_4)\cdot[\widetilde{Y}]}{(1+\mathcal{Y}_2)(1+V)(1+[\widetilde{Y}])}f_4^*c(\widetilde{\mathbb{P}}^{(3)}(\mathscr{E})) \\
                &=&\frac{(1+E_4)(1+\mathcal{Y}_2-E_4)(1+V-E_4)(T-E_4)}{(1+\mathcal{Y}_2)(1+V)(1+T-E_4)}f_4^*c(\widetilde{\mathbb{P}}^{(3)}(\mathscr{E})),
\end{eqnarray*}
\\
where $T=3H+6L-2E_1-2E_2-E_3$ and $\iota:\widetilde{Y}\hookrightarrow \widetilde{\mathbb{P}}(\mathscr{E})$ is the inclusion. By the projection formula,

\[
{f_4}_*c(\widetilde{Y})={f_4}_*\left(\frac{(1+E_4)(1+\mathcal{Y}_2-E_4)(1+V-E_4)(T-E_4)}{(1+T-E_4)}\right)\frac{c(\widetilde{\mathbb{P}}^{(3)}(\mathscr{E}))}{(1+\mathcal{Y}_2)(1+V)}.
\]
\\
Now let $C=\frac{(1+E_4)(1+\mathcal{Y}_2-E_4)(1+V-E_4)(T-E_4)}{(1+T-E_4)}$ and let $\alpha_0, \cdots, \alpha_3$ be the classes obtained by expanding the numerator of $C$ as a polynomial in $E_4$:

\[
(1+E_4)(1+\mathcal{Y}_2-E_4)(1+V-E_4)(T-E_4)=\alpha_0+\alpha_1E_4+\alpha_2E_4^2+\alpha_3E_4^3-E_4^4
\]
\\
Lemma~\ref{BC2} (along with the projection formula and the fact that $f_*f^*\alpha=\alpha$ for any blowup $f$) then gives

\[
{f_4}_*(C)=\alpha_0+\frac{\alpha_2}{V-\mathcal{Y}_2}\left(\frac{\mathcal{Y}_2^2V}{1+T-\mathcal{Y}_2}-\frac{V^2\mathcal{Y}_2}{1+T-V}\right)+
\]
\[
+\frac{\alpha_3}{V-\mathcal{Y}_2}\left(\frac{\mathcal{Y}_2^3V}{1+T-\mathcal{Y}_2}-\frac{V^3\mathcal{Y}_2}{1+T-V}\right)-\frac{1}{V-\mathcal{Y}_2}\left(\frac{\mathcal{Y}_2^4V}{1+T-\mathcal{Y}_2}-\frac{V^4\mathcal{Y}_2}{1+T-V}\right).
\]
\\
Computing three more pushforwards in the same manner then yields the stringy Chern class of $Y$:

\[
c_{str}(Y)=c(Z)+X\cdot [Y_{sing}],
\]
\\
where $\psi:Z\to B$ is a smooth Weierstrass fibration, $[Y_{sing}]=(H+2L)(H+3L)D$ is the class of the singular locus of $Y$ and $X$ is a (lengthy) rational expression in $H$, $L$ and $D:=D_{GUT}$ multiplied by $\pi^*c(B)$ which we do not write explicitly. The stringy Milnor class of $Y$ is then $X\cdot [Y_{sing}]$ and we recall that 

\[
c(Z)=\frac{(1+H)(1+H+2L)(1+H+3L)(3H+6L)}{1+3H+6L}\pi^*c(B).
\]
\\
As the pushforward to the base of $c(Z)$ has been computed in e.g. \cite{AE1}, computing $\varphi_*c_{str}(Y)$ amounts to computing $\varphi_*\mathcal{M}_{str}(Y)=\varphi_*(X\cdot [Y_{sing}])$. For this, we view $X\cdot [Y_{sing}]$ as a class in $A^*\mathbb{P}(\mathscr{E})$ and push it forward to the base via the pushforward formula for classes in a projective bundle first derived in \cite{James}. To apply the pushforward formula of \cite{James} to the case of a projective bundle of the form $\mathbb{P}(\mathscr{O}\oplus \mathscr{L}^2\oplus \mathscr{L}^3)$, we first expand $X\cdot [Y_{sing}]$ as a series in $H$:

\[
X\cdot [Y_{sing}]=\nu_0+\nu_1H+\cdots
\]
\\
Then we consider the following expression (viewed as a function of $H$):

\[
F(H)=\frac{(X\cdot [Y_{sing}]-\gamma)}{H^2},
\]
\\
where $\gamma=\nu_0+\nu_1H$. Then the pushforward to the base of $X\cdot [Y_{sing}]$ is precisely

\[
3\cdot \left.F\right|_{-3L}-2\cdot \left.F\right|_{-2L}=W\cdot c(D),
\]
\\
where $W$ is as given in the conclusion of Proposition~\ref{SVW}. Thus

\[
\varphi_*c_{str}(Y)=\psi_*c(Z)+W\cdot c(D),
\]
\\
from which Proposition~\ref{SVW} immediately follows. Upon integration of $\varphi_*c_{str}(Y)$ we arrive at the following expressions for stringy Euler characterstics for SU(5), SO(10) and $\text{E}_6$ elliptic fibrations in terms of Chern classes of $B$, $D$ and $\mathscr{L}$:

\begin{table}[hbt]
\begin{center}
\begin{tabular}{|c|c|}
\hline
\text{dim}(B) & $\chi_{str}(Y)$ \\
\hline
1 & $12L$  \\
\hline
2 & $12Lc_1-72L^2+80LD-30D^2$  \\
\hline
3 & $12Lc_2-72L^2c_1+432L^3+(80Lc_1-830L^2)D+(555L-30c_1)D^2-120D^3$ \\
\hline
\end{tabular}
\end{center}
\caption{Stringy Euler characteristics of singular \text{SU}(5) elliptic fibrations.}
\end{table}

\begin{table}[hbt]
\begin{center}
\begin{tabular}{|c|c|}
\hline
\text{dim}(B) & $\chi_{str}(Y)$ \\
\hline
1 & $12c_1$  \\
\hline
2 & $80c_1D-60c_1^2-30D^2$  \\
\hline
3 & $288+360c_1^3-750c_1^2D+525c_1D^2-120D^3$ \\
\hline
\end{tabular}
\end{center}
\caption{Stringy Euler characteristics of singular Calabi-Yau \text{SU}(5) elliptic fibrations.}
\end{table}

\begin{table}[hbt]
\begin{center}
\begin{tabular}{|c|c|}
\hline
\text{dim}(B) & $\chi_{str}(Y)$ \\
\hline
1 & $12L$  \\
\hline
2 & $12Lc_1-72L^2+84LD-32D^2$  \\
\hline
3 & $12Lc_2-72L^2c_1+432L^3+(84Lc_1-840L^2)D+(560L-32c_1)D^2-120D^3$ \\
\hline
\end{tabular}
\end{center}
\caption{Stringy Euler characteristics of singular \text{SO}(10) elliptic fibrations.}
\end{table}

\begin{table}[hbt]
\begin{center}
\begin{tabular}{|c|c|}
\hline
\text{dim}(B) & $\chi_{str}(Y)$ \\
\hline
1 & $12c_1$  \\
\hline
2 & $84c_1D-60c_1^2-32D^2$  \\
\hline
3 & $288+360 c_1^3-756c_1^2D+528c_1D^2-120D^3$ \\
\hline
\end{tabular}
\end{center}
\caption{Stringy Euler characteristics of singular Calabi-Yau \text{SO}(10) elliptic fibrations.}
\end{table}

\begin{table}[hbt]
\begin{center}
\begin{tabular}{|c|c|}
\hline
\text{dim}(B) & $\chi_{str}(Y)$ \\
\hline
1 & $12L$  \\
\hline
2 & $12Lc_1-72L^2+90LD-36D^2$  \\
\hline
3 & $12Lc_2-72L^2c_1+432L^3+(-864L^2+90Lc_1)D+(585L-36c_1)D^2-126D^3$ \\
\hline
\end{tabular}
\end{center}
\caption{Stringy Euler characteristics of singular $\text{E}_6$ elliptic fibrations.}
\end{table}

\begin{table}[hbt]
\begin{center}
\begin{tabular}{|c|c|}
\hline
\text{dim}(B) & $\chi_{str}(Y)$ \\
\hline
1 & $12c_1$  \\
\hline
2 & $90c_1D-60c_1^2-36D^2$  \\
\hline
3 & $288+360 c_1^3-774c_1^2D+549c_1D^2-126D^3$ \\
\hline
\end{tabular}
\end{center}
\caption{Stringy Euler characteristics of singular Calabi-Yau $\text{E}_6$ elliptic fibrations.}
\end{table}

\section{Hirzebruch series of a small resolution $\widetilde{Y}$} Let $X$ be a complex smooth variety, then its \emph{Hirzebruch series} is defined as

\[
\mathscr{H}_{y}(X)=\mathscr{H}_0(X)+\mathscr{H}_1(X)y+\mathscr{H}_2(X)y^2+\cdots=\prod_{i=1}^{\text{dim}(X)}(1+ye^{-\lambda_i})\frac{\lambda_i}{1-e^{-\lambda_i}},
\]
\\ 
where $\lambda_i$ are the Chern roots of the tangent bundle of $X$ and $\mathscr{H}_q(X)=\text{ch}(\Omega_{X}^q(X))\text{td}(X)$, i.e., the Chern character of the $q$th exterior power of the cotangent bundle of $X$ multiplied by the Todd class of $X$. Integrating each term in the Hirzebruch series then yields Hirzebruch's $\chi(y)$ characteristic \cite{H}:

\[
\chi(y)=\chi_0+\chi_1y+\chi_2y^2+\cdots,
\]
\\
where $\chi_q=\int_X \mathscr{H}_q(X)$. By the celebrated Hirzebruch-Riemann-Roch theorem (later generalized by Grothendieck), 

\[
\chi_q=\sum_{i=0}^{\text{dim}(X)}(-1)^ih^{i,q},
\]
\\
where $h^{p,q}$ are the Hodge numbers of $X$. As the Hodge numbers of string vacua are particularly important in the context of string theory, we were naturally motivated to compute Hirzebruch series of GUT vacua. Now it turns out that the Hirzebruch series of a variety satisfies many of the same properties as its Chern polynomial, thus many formulas for Chern classes may be immediately converted into a formula for a Hirzebruch series, which we now explain. 

Let $\mathscr{E}\to X$ be a vector bundle over a smooth variety. Define the Chern-ext character of $\mathscr{E}$ to be

\[
\text{ch}_{ext}(\mathscr{E})=1+\text{ch}(\mathscr{E})y+\text{ch}(\Lambda^2\mathscr{E})+\cdots,
\]
\\ 
which can be interpreted as just the Chern character of the \emph{total $\lambda$-class} of $\mathscr{E}$ \cite{FL}. In \cite{FvH} it was shown that if 

\[
0\longrightarrow \mathscr{A} \longrightarrow \mathscr{B} \longrightarrow \mathscr{C} \longrightarrow 0
\]
\\
is an exact sequence of vector bundles, then $\text{ch}_{ext}$ satisfies the \emph{Whitney property}, i.e.,

\[
\text{ch}_{ext}(\mathscr{B})=\text{ch}_{ext}(\mathscr{A})\text{ch}_{ext}(\mathscr{C}).
\]
\\
The fact that the usual Chern character is additive with respect to exact sequences is encoded in the degree one piece of the formula above for the Chern-ext character. As is well known, Todd classes also satisfy the Whitney property. Thus we can associate with any vector bundle $\mathscr{E}\to X$ a \emph{Hirzebruch series}

\[
\mathscr{H}_{y}(\mathscr{E})=\text{ch}_{ext}(\mathscr{E}^{\vee})\text{td}(\mathscr{E}),
\]
\\ 
and this Hirzebruch series satisfies the Whitney property as well, i.e., given an exact sequence of vector bundles $0\longrightarrow \mathscr{A} \longrightarrow \mathscr{B} \longrightarrow \mathscr{C} \longrightarrow 0$ we have

\[
\mathscr{H}_{y}(\mathscr{B})=\mathscr{H}_{y}(\mathscr{A})\mathscr{H}_{y}(\mathscr{C}).
\]
\\
The Hirzebruch series of a smooth variety $X$ is then $\mathscr{H}_{y}(TX)\cap [X]$, and if $X\hookrightarrow M$ is a regular embedding into a smooth variety $M$ then adjunction holds, i.e.,

\[
\mathscr{H}_{y}(X)=\frac{\mathscr{H}_{y}(TM)}{\mathscr{H}_{y}(N)}\cap [X],
\]
\\
where $N$ is the bundle which restricts to the normal bundle of $X$ in $M$. Moreover, since the Chern character and Todd class are both defined in terms of symmetric functions in Chern roots, the Hirzebruch series satisfies all the properties which characterize Chern classes (such as the projection formula and functoriality  $\mathscr{H}_y(f^*\mathscr{E})=f^*\mathscr{H}_y(\mathscr{E})$) but with a different normalization condition. The normalization condition can be recovered from the definition: If $\mathscr{L}\to X$ is a line bundle and $L=c_1(\mathscr{L})$ is its first Chern class then 

\[
\mathscr{H}_{y}(\mathscr{L})=\frac{(1+ye^{-L})L}{1-e^{-L}}.
\]
\\
Thus given a Chern class formula that was derived using only its characterizing properties then you immediately have a formula for its Hirzebruch series as well\footnote{However, one must be careful when converting Chern class formulas into Hirzebruch series formulas, as $c(\mathscr{O})=1$ while $\mathscr{H}_y(\mathscr{O})=1+y$.}. For example, the Hirzebruch series version of Aluffi's formula for blowing up Chern classes is the following

\begin{lemma}Let $f:\widetilde{X}\to X$ be the blowup of a smooth variety $X$ along a smooth complete intersection $V:(F_1=F_2=\cdots =F_k=0)\subset X$, let $E\in A^*\widetilde{X}$ be the class of the exceptional divisor and let $U_i\in A^{*}X$ be the class of $F_i=0$. Then

\[
\mathscr{H}_y(T_{\widetilde{X}})=\frac{\mathscr{H}_y(\mathscr{O}(E))\mathscr{H}_y(\mathscr{O}(f^*U_1-E))\cdots \mathscr{H}_y(\mathscr{O}(f^*U_k-E))}{(1+y)\mathscr{H}_y(\mathscr{O}(f^*U_1))\cdots \mathscr{H}_y(\mathscr{O}(f^*U_k))}\cdot f^*\mathscr{H}_y(T_X).
\]
\end{lemma}

As such, we may compute Hirzebruch series of small resolutions of GUT vacua just as we computed their Chern classes \emph{mutatis mutandis}. Futhermore, if $g:X\to V$ is a proper morphism we define $g_*\mathscr{H}_y(X)$ in the obvious way:

\[
g_*\mathscr{H}_y(X)=g_*\mathscr{H}_0(X)+g_*\mathscr{H}_1(X)y+\cdots
\]
\\
Moreover, if $f:\widetilde{Y}\to Y$ is a small resolution of a GUT elliptic fibration we will then (tentatively) refer to $f_*\mathscr{H}_y(\widetilde{Y})$ as the \emph{stringy Hirzebruch series of} $Y$, which we will denote by $\mathscr{H}_y^{str}(Y)$. The analogue of Proposition~\ref{SVW} is then recorded in the following

\begin{prop}\label{H} Let $\varphi:Y\to B$ be an \emph{SU(5)}, \emph{SO(10)} or $\emph{E}_6$ elliptic fibration with $\text{\emph{dim}}(B)\leq3$, $D=D_{\emph{GUT}}\in A^*B$ and let $L=c_1(\mathscr{L})$. Then 

\[
\varphi_*\mathscr{H}_y^{str}(Y)=\left(1-y+\frac{(1+y)(ye^{-5L}-e^{-L})}{(1+ye^{-6L})}+P\right)\cdot \mathscr{H}_{y}(B), 
\]
\\
where $P=P_1y+P_2y^2+P_3y^3+P_4y^4$ with $P_i$ a polynomial in $D$ and $L$ for $i=1,\ldots,4$, each of which is listed below. 
\end{prop}

\begin{table}[hbt]
\begin{center}
\begin{tabular}{|c|c|c|c|}
\hline
 & SU(5) case \\
\hline
$P_1$ & $20D^3+15D^2-95LD^2-40LD+145L^2D$  \\
\hline
$P_2$ & $-140D^3-45D^2+650LD^2+120DL-975DL^2$   \\
\hline
$P_3$ & $380D^3+75D^2-1760LD^2+2635DL^2-200LD$   \\
\hline
$P_4$ & $-740D^3-105D^2+280LD-5125DL^2+3425LD^2$   \\
\hline
\end{tabular}
\end{center}
\caption{Polynomial coefficients of $P$ in the SU(5) case as in Proposition~\ref{H}.}
\end{table}

\begin{table}[hbt]
\begin{center}
\begin{tabular}{|c|c|}
\hline
 & SO(10) case  \\
\hline
$P_1$ & $20D^3+16D^2-96LD^2+147L^2D-42LD$  \\
\hline
$P_2$ & $-140D^3-48D^2+656LD^2-987L^2D+126LD$   \\
\hline
$P_3$ & $380D^3+80D^2-1776LD^2+2667L^2D-210LD$   \\
\hline
$P_4$ & $-740D^3-112D^2+294LD-5187L^2D+3456LD^2$   \\
\hline
\end{tabular}
\end{center}
\caption{Polynomial coefficients of $P$ in the SO(10) case as in Proposition~\ref{H}.}
\end{table}

\begin{table}[hbt]
\begin{center}
\begin{tabular}{|c|c|c|c|}
\hline
 & $\text{E}_6$ case \\
\hline
$P_1$ & $21D^3+18D^2+\frac{303}{2}L^2D-45LD-\frac{201}{2}LD^2$ \\
\hline
$P_2$ & $-147D^3-54D^2-\frac{2031}{2}L^2D+135LD+\frac{1371}{2}LD^2$  \\
\hline
$P_3$ & $399D^3+90D^2-\frac{3711}{2}LD^2+\frac{5487}{2}L^2D-225LD$   \\
\hline
$P_4$ & $-777D^3-126D^2+\frac{7221}{2}LD^2-\frac{10671}{2}L^2D+315LD$   \\
\hline
\end{tabular}
\end{center}
\caption{Polynomial coefficients of $P$ in the $\text{E}_6$ case as in Proposition~\ref{H}.}
\end{table}

\section{A cohomological speculation} To our knowledge no one (including ourselves) as of yet has computed the Hodge numbers of small resolutions of GUT vacua (which we recall coincide with the stringy Hodge numbers of singular GUT elliptic fibrations). By a theorem of Kontsevich \cite{Kontsevich}, birationally equivalent Calabi-Yau varieties have the same Hodge numbers, thus the Hodge numbers of any crepant resolution of a singular GUT elliptic fibration will coincide. It would be nice if one could use the Lefschetz hyperplane theorem to relate the upper cohomology of the total space of a smooth elliptic fibration which is embedded in a (possibly blown up) projective bundle to the cohomology of the ambient projective bundle in which it resides, but unfortunately the total space of such an elliptic fibration is not an ample divisor, thus the hypotheses of the Lefschetz theorem are not satisfied. However, in cases where the Hodge numbers of elliptic fibrations are known (e.g. when the base of the fibration is a Fano toric variety) it is indeed the case that the upper cohomology of the total space of the elliptic fibration $Y$ coincides with the ambient projective bundle in which it resides. We speculate that this is due to the fact that even though $Y$ is not an ample divisor, perhaps

\[
\tag{\ddag}
H^q(\mathbb{P}(\mathscr{E}),\Omega_{\mathbb{P}(\mathscr{E})}^p(-Y))=H^q(Y,\Omega_Y^{p-1}(-Y))=0
\]
\\
for $p+q<\text{dim}(\mathbb{P}(\mathscr{E}))$, which would then imply the result of Lefschetz\footnote{In the algebraic proof of the Lefschetz hyperplane theorem, it is property $(\ddag)$ which ultimately yields the conclusion of the theorem. Though the ampleness (or `amplitude') of $Y$ is sufficient (by Kodaira's vanishing theorem) to conclude property $(\ddag)$, it is not necessary.}. Moreover, we speculate that property $(\ddag)$ holds for $Y$ \emph{any} hypersurface in a (possibly blown up) projective bundle which restricts to a hypersurface (or local complete intersection) on each fiber. Assuming this is the case, then the Hodge numbers for any smooth Weierstrass fibration or resolution of a singular Weierstrass fibration would be easily computable in terms of data on the base using Hodge-Deligne polynomials and Hirzebruch-Riemann-Roch\footnote{The same would then follow for the non-Weierstrass fibrations referred to as $E_7$, $E_6$ and $D_5$ fibrations as well \cite{AE1}\cite{AE2}\cite{EF}.}. For example, let $\varphi:Y\to B$ be a smooth Calabi-Yau Weierstrass fibration over a three dimensional base as defined in \S\ref{intro}. Then by Hirzebruch-Riemann-Roch and the pushforward formula of \cite{James} we have the following relations between the Hodge numbers of $Y$ \cite{EF}:

\[
\tag{\ddag\ddag}
h^{1,2}-h^{1,1}-h^{1,3}=-40-60c_1(B)^3, \quad h^{2,2}-2h^{2,1}=204+240c_1(B)^3
\]
\\
Now if property $(\ddag)$ indeed holds, then $h^{1,1}$ and $h^{1,2}$ coincide with the corresponding Hodge numbers of the total space of the projective bundle $\pi:\mathbb{P}(\mathscr{E})\to B$ in which it resides, which are easily computed in terms of the Hodge numbers of $B$ using Hodge-Deligne polynomials. The Hodge numbers $h^{2,2}$ and $h^{1,3}$ then could immediately be written in terms of the Hodge numbers of $B$ and $c_1(B)^3$ by equations $(\ddag\ddag)$.

To be more precise, we recall that if $X$ a smooth projective variety then the Hodge-Deligne polynomial of $X$ is simply

\[
E(X)=E(X;u,v):=\sum_{p,q}(-1)^{p+q}h^{p,q}(X)u^pv^q.
\]
\\
For example we have

\[
E(\mathbb{P}^n)=1+uv+\cdots+(uv)^n.
\]
\\
It is a fact that for a Zariski locally trivial fibration $f:X\to Z$ with fiber $F$ we have $E(X)=E(Z)\cdot E(F)$, thus the Hodge-Deligne polynomial of the ambient projective bundle $\mathbb{P}(\mathscr{E})$ is

\[
E(\mathbb{P}(\mathscr{E}))=(1+uv+(uv)^2)\cdot E(B).
\]
\\
In particular if $B$ is Fano\footnote{The Fano assumption is not necessary, though (almost) Fano varieties are natural candidates for F-theory bases.}, then

\[
E(B)=1+h^{1,1}(B)uv-h^{1,2}(B)(u^2v+uv^2)+h^{1,1}(B)u^2v^2+u^3v^3,
\]
\\
thus

\[
h^{1,1}(\mathbb{P}(\mathscr{E}))=h^{1,1}(B)+1, \quad h^{1,2}(\mathbb{P}(\mathscr{E}))=h^{1,2}(B).
\]
\\
Assuming property $(\ddag)$ holds, then equations $(\ddag\ddag)$ yield the following formulas for the non-trivial Hodge numbers of a smooth Weierstrass fibration $Y$ solely in terms of the Hodge numbers of $B$ and $c_1(B)^3$:

\[
\begin{cases}
h^{1,1}(Y)=h^{1,1}(B)+1 \\
h^{1,2}(Y)=h^{1,2}(B) \\
h^{1,3}(Y)=39+60c_1(B)^3-h^{1,1}(B)+h^{1,2}(B) \\
h^{2,2}(Y)=204+240c_1(B)^3+2h^{1,2}(B)
\end{cases}
\]
\\
These formulas are correct for $B$ any of the 18 Fano toric varieties \cite{KLRY}. We speculate that not only are these formulas correct for any smooth Calabi-Yau Weierstrass fibration over a Fano base, but similar formulas may be derived analogously for small resolutions of singular GUT elliptic fibrations, as one can easily compute Hodge-Deligne polynomials of blowups provided the Hodge-Deligne polynomial of the variety which is being blown up and the Hodge-deligne polynomial of the center of the blowup are known.

\appendix

\section{Small resolutions of GUT vacua}\label{A1}
We now recall the resolution procedures for SU(5), SO(10) and $\text{E}_6$ elliptic fibrations as first presented in the physics literature \cite{EY}\cite{SO(10)}\cite{E6} (we follow very closely the global description of a resolution presented in \cite{MS}). 
\\
\subsection{SU(5) resolution(s)}
\noindent\emph{First blowup}: Let $f_1:\widetilde{\mathbb{P}}^{(1)}(\mathscr{E})\to \mathbb{P}(\mathscr{E})$ be the blowup of $\mathbb{P}(\mathscr{E})$ along the singular locus $Y_{sing}:(x=y=w=0)$ of $Y$ and denote the exceptional divisor by $E_1$. The sections $x$, $y$ and $w$ then pullback as

\[
x\mapsto \delta_1x_1, \quad y\mapsto \delta_1y_1, \quad w\mapsto \delta_1w_1,
\]
\\
where $\delta_1$ is a regular section $\mathscr{O}(E_1)$ (this also determines the classes $[x_1=0]$, $[y_1=0]$ and $[w_1=0]$). The class of the proper transform $Y^{(1)}$ of $Y$ is then $[Y^{(1)}]=3H+6L-2E_1$, where $H=c_1(\mathscr{O}_{\mathbb{P}(\mathscr{E})}(1))$ and $L=c_1(\mathscr{L})$. When $L=c_1(B)$, the proper transform $Y^{(1)}$ of $Y$ is an anti-canonical divisor of $\widetilde{\mathbb{P}}^{(1)}(\mathscr{E})$, as $c_1(\widetilde{\mathbb{P}}^{(1)}(\mathscr{E}))=f_1^*c_1(\mathbb{P}(\mathscr{E}))+(1-3)E_1=3H+c_1(B)+5L-2E_1$. By Lemma~\ref{BC} $(\dag)$, we have

\[
c(\widetilde{\mathbb{P}}^{(1)}(\mathscr{E}))=\frac{(1+E_1)(1+H+2L-E_1)(1+H+3L-E_1)(1+D-E_1)}{(1+H+2L)(1+H+3L)(1+D)}f_1^*c(\mathbb{P}(\mathscr{E})) \]
\\
where $D:=c_1(\mathscr{L}_{\text{GUT}})$.

\vspace{0.22in}

\noindent\emph{Second blowup}: Let $f_2:\widetilde{\mathbb{P}}^{(2)}(\mathscr{E})\to \widetilde{\mathbb{P}}^{(1)}(\mathscr{E})$ be the blowup of $\widetilde{\mathbb{P}}^{(1)}(\mathscr{E})$ along

\[
W^{(2)}:(x_1=y_1=\delta_1=0)\subset \widetilde{\mathbb{P}}^{(1)}(\mathscr{E}).
\]
\\
Denote the classes $[x_1=0]=H+2L-E_1$ and $[y_1=0]=H+3L-E_1$ by $\mathcal{X}_1$ and $\mathcal{Y}_1$ respectively, and denote the class of the exceptional divisor by $E_2$. The sections $x_1$, $y_1$ and $\delta_1$ then pullback as

\[
x_1\mapsto \delta_2x_{2}, \quad y_1\mapsto \delta_2y_2, \quad \delta_1\mapsto \delta_2\zeta_2,
\]
\\
where $\delta_2$ is a regular section of $\mathscr{O}(E_2)$. The inverse image of $Y$ through the first two blowups takes the form

\[
(f_1\circ f_2)^{-1}(Y): \left(\zeta_2^2\delta_2^4\cdot(y_2\tau z-\zeta_2\delta_2\chi)=0\right) \subset \widetilde{\mathbb{P}}^{(2)}(\mathscr{E}),
\]
\\
where $\tau =y_2+\zeta_2\beta_3w_1^2z+a_1x_2$ and $\chi=x_2^3\delta_2+\beta_4w_1x_2^2z+\zeta_2\beta_2w_1^3x_2z^2+\zeta_2^2\beta_0w_1^5z^3$. The proper transform $Y^{(2)}$ is then given by $\{y_2\tau z-\zeta_2\delta_2\chi=0\}$ and $[Y^{(2)}]=[Y^{(1)}]-2E_2=3H+6L-2E_1-2E_2$, which again is an anti-canonical divisor of $\widetilde{\mathbb{P}}^{(2)}(\mathscr{E})$ when $L=c_1(B)$. By Lemma~\ref{BC} $(\dag)$, we have

\[
c(\widetilde{\mathbb{P}}^{(2)}(\mathscr{E}))=\frac{(1+E_2)(1+\mathcal{X}_2)(1+\mathcal{Y}_2)(1+E_1-E_2)}{(1+\mathcal{X}_1)(1+\mathcal{Y}_1)(1+E_1)}f_2^*c(\widetilde{\mathbb{P}}^{(1)}(\mathscr{E})), \]
\\
where $\mathcal{X}_2=\mathcal{X}_1-E_2$ and $\mathcal{Y}_2=\mathcal{Y}_1-E_2$.
\vspace{0.22in}

\noindent\emph{Third and fourth blowups}: The third and fourth blowups are obtained by blowing up codimension two loci which results in a small resolution of $Y$. Let $f_3:\widetilde{\mathbb{P}}^{(3)}(\mathscr{E})\to \widetilde{\mathbb{P}}^{(2)}(\mathscr{E})$ be the blowup of $\widetilde{\mathbb{P}}^{(2)}(\mathscr{E})$ along

\[
W^{(3)}:(y_2=u=0)\subset \widetilde{\mathbb{P}}^{(2)}(\mathscr{E}),
\]
\\
and let $f_4:\widetilde{\mathbb{P}}(\mathscr{E})\to \widetilde{\mathbb{P}}^{(3)}(\mathscr{E})$ be the blowup of $\widetilde{\mathbb{P}}^{(3)}(\mathscr{E})$ along

\[
W^{(4)}:(\tau=v=0)\subset \widetilde{\mathbb{P}}^{(3)}(\mathscr{E}).
\]
\\
We denote the exceptional divisors by $E_3$ and $E_4$ respectively. The six resolutions are then obtained by taking the proper transform of $Y^{(2)}$ under these blowups corresponding to the following choices for the values of $u$ and $v$:

\begin{table}[hbt]
\begin{center}
\begin{tabular}{|c|c|c|c|c|}
\hline
                           & $u$ & $v$ & $[u=0]$ & $[v=0]$  \\
\hline
$1^{\text{st}}$ resolution & $\delta_2$ & $\zeta_2$ & $E_2$ & $E_1-E_2$ \\
\hline
$2^{\text{nd}}$ resolution & $\delta_2$ & $\chi$ & $E_2$ & $3\mathcal{X}_2+E_2$ \\
\hline
$3^{\text{rd}}$ resolution & $\zeta_2$ & $\delta_2$ & $E_1-E_2$ & $E_2$ \\
\hline
$4^{\text{th}}$ resolution & $\zeta_2$ & $\chi$ & $E_1-E_2$ & $3\mathcal{X}_2+E_2$ \\
\hline
$5^{\text{th}}$ resolution & $\chi$ & $\delta_2$ & $3\mathcal{X}_2+E_2$ & $E_2$ \\
\hline
$6^{\text{th}}$ resolution & $\chi$ & $\zeta_2$ & $3\mathcal{X}_2+E_2$ & $E_1-E_2$ \\
\hline

\end{tabular}
\end{center}
\caption{Choices for $u$ and $v$ (and their classes) corresponding to each of the six small resolutions of SU(5) vacua.}
\end{table}

The sections $y_2$, $u$, $\tau$, and $v$  then pullback as

\[
y_2\mapsto \delta_3y_3, \quad u\mapsto \delta_3u_3, \quad \tau\mapsto \delta_4\tilde{y}, \quad v\mapsto \delta_4v_4,
\]
\\
where $\delta_3$ and $\delta_4$ are regular sections of $\mathscr{O}(E_3)$ and $\mathscr{O}(E_4)$ respectively. By Lemma~\ref{BC} $(\dag)$, we have

\[
c(\widetilde{\mathbb{P}}^{(3)}(\mathscr{E}))=\frac{(1+E_3)(1+\mathcal{Y}_2-E_3)(1+U-E_3)}{(1+\mathcal{Y}_2)(1+U)}f_3^*c(\widetilde{\mathbb{P}}^{(2)}(\mathscr{E})),
\]
\\
and

\[
c(\widetilde{\mathbb{P}}(\mathscr{E}))=\frac{(1+E_4)(1+\mathcal{Y}_2-E_4)(1+V-E_4)}{(1+\mathcal{Y}_2)(1+V)}f_4^*c(\widetilde{\mathbb{P}}^{(3)}(\mathscr{E})),
\]
\\
where $U$ and $V$ are the classes of $\{u=0\}$ and $\{v=0\}$ respectively. The pullback of the defining equation for $Y$ under the four blowups (corresponding to the first resolution in Table 1) is

\[
\delta_4^3v_4^2\delta_3^5u_3^4(y_4y_3z-v_4u_3\chi)=0,
\]
\\
thus the proper transform of $Y$ is given by

\[
\widetilde{Y}_{\text{SU(5)}}:\left(y_4y_3z-v_4u_3\chi=0\right) \subset \widetilde{\mathbb{P}}(\mathscr{E}),
\]
\\
and is of class
\[
[\widetilde{Y}_{\text{SU(5)}}]=[Y^{(3)}]-E_4=[Y^{(2)}]-E_3-E_4=3H+6L-2E_1-2E_2-E_3-E_4.
\]
\\
Moreover, $\widetilde{Y}_{\text{SU(5)}}$ is easily checked to be a smooth hypersurface of $\widetilde{\mathbb{P}}(\mathscr{E})$ which is $K$-equivalent to the original singular fourfold $Y$, and is an anti-canonical divisor when $L=c_1(B)$.

\subsection {SO(10) resolution} For the SO(10) case we follow the procedure first presented in \cite{SO(10)} which requires five blowups for the small resolution of $Y$, and only give the pullback maps associated with each blowup as everything else follows exactly as in the SU(5) case \emph{mutatis mutandis}.
\\

\noindent\emph{First and second blowups}: The first and second blowups are precisely the same as in the SU(5) case and we use the same notations.
\\

\noindent\emph{Third blowup}: Let $g_3:\widetilde{\mathbb{P}}^{(III)}(\mathscr{E})\to \widetilde{\mathbb{P}}^{(2)}(\mathscr{E})$ be the blowup of $\widetilde{\mathbb{P}}^{(2)}(\mathscr{E})$ along

\[
W'^{(3)}:(y_{2}=\zeta_2=\delta_2=0)\subset \widetilde{\mathbb{P}}^{(II)}(\mathscr{E}),
\]
\\
and denote the exceptional divisor by $E_3$. The sections $y_2$, $\zeta_2$ and $\delta_2$ then pullback as

\[
y_2\mapsto \delta_3y_3, \quad \zeta_2\mapsto \delta_3\zeta_3, \quad \delta_2\mapsto \delta_3\xi_3,
\]
\\
where $\delta_3$ is a regular section of $\mathscr{O}(E_3)$. The inverse image of the SO(10) fibration through the first three blowups takes the form

\[
(g_1\circ g_2\circ g_3)^{-1}(Y): \left(\delta_3^8\xi_3^4\zeta_3^2\cdot(y_3\tau_3z-\zeta_3\xi_3\chi_3)=0\right) \subset \widetilde{\mathbb{P}}^{(III)}(\mathscr{E}),
\]
\\
where $\tau_3=\zeta_3\vartheta w_1x_2\xi_3\delta_3+y_3+\zeta_3\beta_3w_1^2z$ and $\chi_3=\zeta_3^2\beta_0w_1^5z^3\delta_3^2+\delta_3x_2^3\xi_3+\delta_3\zeta_3\beta_2w_1^3x_2z^2+\beta_4w_1x_2^2z$. The proper transform $\widetilde{Y}^{(3)}$ is then given by $\{y_3\tau_3z-\zeta_3\xi_3\chi_3=0\}$. 
\\

\noindent\emph{Fourth and fifth blowups}: The fourth and fifth blowups are obtained by blowing up codimension two loci which results in a small resolution of the SO(10) fibration. Let $g_4:\widetilde{\mathbb{P}}^{(IV)}(\mathscr{E})\to \widetilde{\mathbb{P}}^{(III)}(\mathscr{E})$ be the blowup of $\widetilde{\mathbb{P}}^{(2)}(\mathscr{E})$ along

\[
W'^{(4)}:(y_3=\zeta_3=0)\subset \widetilde{\mathbb{P}}^{(III)}(\mathscr{E}),
\]
\\
and let $g_5:\widetilde{\mathbb{P}^*}(\mathscr{E})\to \widetilde{\mathbb{P}}^{(IV)}(\mathscr{E})$ be the blowup of $\widetilde{\mathbb{P}}^{(IV)}(\mathscr{E})$ along

\[
W'^{(5)}:(y_4=\xi_3=0)\subset \widetilde{\mathbb{P}}^{(IV)}(\mathscr{E}),
\]
\\
where $y_4$ is the pullback of $y_3$ minus the exceptional divisor. We denote the exceptional divisors by $E_4$ and $E_5$ respectively. The sections $y_3$, $\zeta_3$, $y_4$ and $\xi_3$  then pullback as

\[
y_3\mapsto \delta_4y_4, \quad \zeta_3\mapsto \delta_4 \zeta_4, \quad y_4\mapsto \delta_5y_5, \quad \xi_3\mapsto \delta_5\xi_5,
\]
\\
where $\delta_4$ and $\delta_5$ are regular sections of $\mathscr{O}(E_4)$ and $\mathscr{O}(E_5)$ respectively. The pullback of the defining equation for the SO(10) fibration under the five blowups is

\[
\delta_4^3\zeta_4^2\delta_3^8\delta_5^5\xi_5^4\cdot(y_5y_4z-\zeta_4\xi_5\chi_3)=0,
\]
thus the proper transform of the SO(10) fibration is then given by

\[
\widetilde{Y}_{\text{SO(10)}}:\left(y_5y_4z-\zeta_4\xi_5\chi_3=0\right) \subset \widetilde{\mathbb{P}'}(\mathscr{E}),
\]
\\
and is of class
\[
[\widetilde{Y}_{\text{SO(10)}}]=3H+6L-2E_1-2E_2-2E_3-E_4-E_5.
\]
\\
Moreover, $\widetilde{Y}_{\text{SO(10)}}$ is easily checked to be a smooth hypersurface of $\widetilde{\mathbb{P}'}(\mathscr{E})$ which is $K$-equivalent to the original singular SO(10) fibration, and is an anti-canonical divisor when $L=c_1(B)$.

\subsection{$\text{E}_6$ resolution} For the $\text{E}_6$ case we follow the procedure first presented in \cite{E6} which requires seven blowups for the small resolution of $Y$, and only give the pullback maps associated with each blowup as everything else follows exactly as in the SU(5) and SO(10) cases \emph{mutatis mutandis}.
\\

\noindent\emph{First three blowups}: The first three blowups are precisely the same as in the SO(10) case and we use the same notations.
\\

\noindent\emph{Fourth blowup}:  Let $h_4:\widetilde{\mathbb{P}}^{(\diamondsuit)}(\mathscr{E})\to \widetilde{\mathbb{P}}^{(III)}(\mathscr{E})$ be the blowup of $\widetilde{\mathbb{P}}^{(III)}(\mathscr{E})$ along

\[
\mathscr{W}^{(4)}:(y_{3}=\zeta_3=\delta_3=0)\subset \widetilde{\mathbb{P}}^{(III)}(\mathscr{E}),
\]
\\
and denote the exceptional divisor by $E_4$. The sections $y_3$, $\zeta_3$ and $\delta_3$ then pullback as

\[
y_3\mapsto \delta_4y_4, \quad \zeta_3\mapsto \delta_4\zeta_4, \quad \delta_3\mapsto \delta_4\xi_4,
\]
\\
where $\delta_4$ is a regular section of $\mathscr{O}(E_4)$. 
\\

\noindent\emph{Final three blowups}:  The final three blowups are obtained by blowing up codimension two loci which results in a small resolution of the $\text{E}_6$ fibration. Let $h_5:\widetilde{\mathbb{P}}^{(\clubsuit)}(\mathscr{E})\to \widetilde{\mathbb{P}}^{(\diamondsuit)}(\mathscr{E})$ be the blowup of $\widetilde{\mathbb{P}}^{(\diamondsuit)}(\mathscr{E})$ along

\[
\mathscr{W}^{(5)}:(y_{4}=\xi_3=0)\subset \widetilde{\mathbb{P}}^{(\diamondsuit)}(\mathscr{E}),
\]
\\
let $h_6:\widetilde{\mathbb{P}}^{(\heartsuit)}(\mathscr{E})\to \widetilde{\mathbb{P}}^{(\clubsuit)}(\mathscr{E})$ be the blowup of $\widetilde{\mathbb{P}}^{(\clubsuit)}(\mathscr{E})$ along

\[
\mathscr{W}^{(6)}:(y_{5}=\xi_4=0)\subset \widetilde{\mathbb{P}}^{(\clubsuit)}(\mathscr{E}),
\]
\\
and let $h_7:\widetilde{\mathbb{P}}^{(\spadesuit)}(\mathscr{E})\to \widetilde{\mathbb{P}}^{(\heartsuit)}(\mathscr{E})$ be the blowup of $\widetilde{\mathbb{P}}^{(\heartsuit)}(\mathscr{E})$ along

\[
\mathscr{W}^{(7)}:(y_{6}=\zeta_4=0)\subset \widetilde{\mathbb{P}}^{(\heartsuit)}(\mathscr{E}),
\]
\\
where $y_i$ denotes the pullback of $y_{i-1}$ minus the exceptional divisor of the corresponding blowup. We denote the exceptional divisors by $E_5$, $E_6$ and $E_7$ respectively. The sections $y_4$, $\xi_3$, $y_5$, $\xi_5$, $y_6$ and $\zeta_4$  then pullback as

\[
y_4\mapsto \delta_5y_5, \quad \xi_3\mapsto \delta_5 \xi_5, \quad y_5\mapsto \delta_6y_6, \quad \xi_4\mapsto \delta_6\xi_6, \quad y_6\mapsto \delta_7y_7, \quad \zeta_4\mapsto \delta_7\zeta_7,
\]
\\
where $\delta_5$, $\delta_6$ and $\delta_7$ are regular sections of $\mathscr{O}(E_5)$, $\mathscr{O}(E_6)$ and $\mathscr{O}(E_7)$ respectively. The pullback of the defining equation for the $\text{E}_6$ fibration under the seven blowups is

\[
\delta_4^{12} \delta_6^8 \xi_6^8 \delta_5^5 \xi_5^4 \delta_7^3 \zeta_7^2\cdot (\zeta_7\xi_5\xi_6\delta_6\tau_7-y_7\delta_5\delta_7z\chi_7)=0,
\]
where $\tau_7=\beta_0w_1^5z^3\delta_4^3\delta_7^2\delta_6\xi_6\zeta_7^2+\zeta_7\delta_7\delta_5\eta w_1^2x_2^2z\delta_4^2\delta_6\xi_6\xi_5+\zeta_7\delta_4\delta_7\beta_2w_1^3x_2z^2+\delta_5x_2^3\xi_5$ and $\chi_7=y_7\delta_5^2+\delta_4\delta_6\xi_6\xi_5\zeta_7\vartheta w_1x_2\delta_5+\zeta_7\beta_3w_1^2z$. The proper transform of the $\text{E}_6$ fibration is then given by

\[
\widetilde{Y}_{\text{E}_6}:\left(\zeta_7\xi_5\xi_6\delta_6\tau_7-y_7\delta_5\delta_7z\chi_7=0\right) \subset \widetilde{\mathbb{P}}^{(\spadesuit)}(\mathscr{E}),
\]
\\
and is of class
\[
[\widetilde{Y}_{\text{E}_6}]=3H+6L-2E_1-2E_2-2E_3-2E_4-E_5-E_6-E_7.
\]
\\
Moreover, $\widetilde{Y}_{\text{E}_6}$ is easily checked to be a smooth hypersurface of $\widetilde{\mathbb{P}}^{(\spadesuit)}(\mathscr{E})$ which is $K$-equivalent to the original singular $\text{E}_6$ fibration, and is an anti-canonical divisor when $L=c_1(B)$.


\thebibliography{99}

\bibitem{Vafa}
C.~Vafa,
``Geometry of Grand Unification'',
arXiv:0911.3008 [math-ph].

\bibitem{BIKMSV}
M.~Bershadsky, K. A.~Intriligator, S.~Kachru, D.R.~Morrison, V.~Sadov and C.~Vafa, 
``Geometric singularities and enhanced gauge symmetries'', Nucl. Phys. B 481, 215 (1996) [arXiv:hep-th/9605200].

\bibitem{KatzVafa}
S.H.~Katz and C.~Vafa, 
``Matter from geometry'', Nucl. Phys. B 497, 146 (1997) [arXiv:hep-th/9606086].

\bibitem{Tate}
J.~Tate, 
``Algorithm for Determining the Type of a Singular Fiber in an Elliptic Pencil'', Modular Functions of One Variable IV, Lecture Notes in Math 476, Springer-Verlag, 1975, 33-52.

\bibitem{CandelasFont}
P.~Candelas and A.~Font, ``Duality between the webs of heterotic and type II vacua'',
Nucl.Phys. B511 (1998) 295Ð325, hep-th/9603170.

\bibitem{BGJW}
R.~Blumenhagen, T. W.~Grimm, B.~Jurke, and T.~Weigand, ``Global F-theory GUTs'', Nucl.
Phys. B829 (2010) 325Ð369, 0908.1784.

\bibitem{EY}
  M.~Esole, S.~-T.~Yau,
  ``Small resolutions of SU(5)-models in F-theory'',
    [arXiv:1107.0733 [hep-th]].
    
\bibitem{SO(10)}
R.~Tatar and W.~Walters,
``GUT theories from Calabi-Yau 4-folds with SO(10) singularities'',  [arXiv:1206.5090 [hep-th]].

\bibitem{E6}
M.~K\"{u}ntzler and S.~Sch\"{a}fer-Nameki,
``G-flux and Spectral Divisors'',
[arXiv:1205.5688 [hep-th]].

\bibitem{Kontsevich}
M.~Kontsevich, Lecture at Universite\'{} Paris-Sud, Orsay, December 7, 1995.

\bibitem{Batyrev}
V.~Batyrev, ``Stringy Hodge numbers of varieties with Gorenstein canonical singularities'', Proc. Taniguchi Symposium 1997, in ``Integrable Systems and Algebraic Geometry'' Kobe/Kyoto 1997, World Sci. Publ. (1999), 1-32.

\bibitem{Aluffi}
P.~Aluffi, ``Modification systems and integration in their Chow groups'', Selecta Math. (N.S.) 11 (2005), no. 2, 155Ð202.

\bibitem{deFernex}
T.~de Fernex, E.~Lupercio, T.~Nevins and B. Uribe, ``Stringy Chern classes of singular varieties'', Adv. Math. 208 (2007), 597-621.

\bibitem{Batyrev2}
V.~Batyrev, ``Non-Archimedean integrals and stringy Euler numbers of log-terminal pairs'', J. Eur. Math. Soc. (JEMS), 1(1):5Ð 33, 1999.

\bibitem{Aluffi2}
P.~Aluffi, ``Weighted Chern-Mather classes and Milnor classes of hypersurfaces'', Singularities and Arrangements, Sapporo-Tokyo 1998, Advanced Studies in Pure Mathematics 29 (2000) 1--20.

\bibitem{AE1}  P.~Aluffi, M.~Esole,
  ``Chern class identities from tadpole matching in type IIB and F-theory'', JHEP {\bf 0903}, 032 (2009)
  [arXiv:0710.2544 [hep-th]].
  
\bibitem{MS}
J.~Marsano, S.~Schafer-Nameki, ``Yukawas, G-flux, and Spectral Covers from Resolved Calabi-Yau's'',  [arXiv:1108.1794 [hep-th]].

\bibitem{James}
  J.~Fullwood,
  ``On generalized Sethi-Vafa-Witten formulas,'' J. Math. Phys. 52, 082304 (2011).
  
\bibitem{Aluffi3}
P.~Aluffi, ``Chern classes of blow-ups'', Math. Proc. Cambridge Philos. Soc. (2008).

\bibitem{Fulton}
W.~Fulton, ``Intersection theory'', Springer-Verlag, Berlin, 1984.

\bibitem{AE2}
P.~Aluffi, M.~Esole, ``New Orientifold Weak Coupling Limits in F-theory'', JHEP 02 (2010) 020.

\bibitem{H}
F.~Hirzebruch, ``Topological methods in algebraic geometry'',
Springer, New York (1966).

\bibitem{FL}
W.~Fulton, S.~Lang, ``Riemann-Roch algebra'', volume 277 of Grundlehren der Mathematischen Wis-
senschaften [Fundamental Principles of Mathematical Sciences]. Springer-Verlag, New York, 1985.

\bibitem{FvH}
J.~Fullwood, M.~van Hoeij, ``On Hirzebruch invariants of elliptic fibrations'', String Math 2011, [Proceedings of Symposia in Pure Mathematics] AMS, 2012.

\bibitem{EF}
M.~Esole, J.~Fullwood, S.T.~Yau, ``D5 elliptic fibrations: Non-Kodaira fibers and new orientifold
limits of F-theory'', [arXiv:1110.6177 [hep-th]].

\bibitem{KLRY}
A.~Klemm, B.~Lian, S. S.~Roan, S. -T.~Yau, ``Calabi-Yau fourfolds for M- and F-theory
compactificationsÓ, Nucl. Phys. B518, 515-574 (1998). [hep-th/9701023].
\end{document}